\documentclass[11pt]{article}  
\usepackage{times}
\usepackage{graphicx}

\usepackage{cite}
\usepackage{url}
\usepackage{mathfont}

\mathcode`O="724F

\newtheorem{theorem}{Theorem}
\newtheorem{lemma}{Lemma}

\DeclareSymbolFont{lasy}{U}{lasy}{m}{n}
\let\Box\undefined
\DeclareMathSymbol\Box{0}{lasy}{"32}
\newcommand{\qed}{\hfill$\Box$\medbreak}
\newenvironment{proof}{\noindent{\bf Proof: }}{\qed}

\let\thickness\theta
\def\geomthickness{\bar\theta}
\def\bookthickness{{\rm bt}}

\begin{document}

\title{Separating Geometric Thickness from Book Thickness}
\author{David Eppstein\thanks{Univ. of California, Irvine, Dept. Inf. \&
Comp. Sci., Irvine, CA 92697. Work supported in part by NSF grant
CCR-9912338.}}
\date{ }
\maketitle

\begin{abstract}
We show that geometric thickness and book thickness are not
asymptotically equivalent: for every $t$, there exists a graph with
geometric thickness two and book thickness $\ge t$.
\end{abstract}

\section{Introduction}

Graph drawing~\cite{DiBEadTam-99} concerns itself with the
geometric layout of abstract graphs, with applications including
information visualization and VLSI design.  The graphs arising from these
applications are frequently impossible to arrange in the plane without
edge crossings; one way of dealing with this problem is to partition the
edges of the graph into a number of planar {\em layers}, each of which
might be drawn using a different color or placed within a different
physical layer of a VLSI circuit.  This leads naturally to the concept of
{\em thickness}, of which there are several variants:

\begin{itemize}
\item The {\em thickness} of a graph $G$, denoted $\thickness(G)$, is the
minimum number of planar subgraphs into which the edges of $G$ can be
partitioned.  Equivalently, it is the minimum number of layers in a
planar drawing of $G$, such that each edge belongs to a single layer, no
two edges in the same layer cross, and edges are allowed to be drawn as
arbitrary curves~\cite{Kai-AMSUH-73}.

\item The {\em geometric thickness} of a graph $G$, denoted
$\geomthickness(G)$, is the minimum number of layers in a planar drawing
of $G$, such that each edge belongs to a single layer, no two edges in
the same layer cross, and edges must be drawn as straight line
segments.  This parameter was introduced under the name ``real linear
thickness'' by Kainen~\cite{Kai-AMSUH-73}, and further studied by
Dillencourt et al.~\cite{DilEppHir-JGAA-00}.

\item The {\em book thickness} of a graph $G$, denoted
$\bookthickness(G)$, can be defined as the minimum number of layers in a planar drawing
of $G$, such that each edge belongs to a single layer, no two edges in
the same layer cross, and edges must be drawn as straight line segments,
with the further restriction that the vertices of $G$ must be placed in
convex position~\cite{BerKai-JCTB-79}.
\end{itemize}

\noindent
It is not difficult to define further variants;
for instance, Wood~\cite{Woo-COMB-01} considers layouts in which each
edge is drawn with at most one bend, at which it may change layers.  For
more results on thickness, see the survey of Mutzel et
al.~\cite{MutOdeSch-GC-98}.

Clearly, from these definitions,
$\thickness(G)\le\geomthickness(G)\le\bookthickness(G)$, and these
inequalities have been shown to be strict for certain
graphs~\cite{DilEppHir-JGAA-00}.  However, it was not known whether
the geometric thickness is asymptotically equivalent to either of the
other two parameters; that is, whether
$\bookthickness(G)=O(\geomthickness(G))$ or whether
$\geomthickness(G)=O(\thickness(G))$.
In this paper, we answer the first of these two questions in the
negative, by exhibiting a family of graphs for which $\geomthickness=2$
but for which $\bookthickness=\omega(1)$.

Our construction is closely
related to the layout method of Wood~\cite{Woo-COMB-01} however rather
than allowing bends in the edges of our drawings we build them into the
input by subdividing the edges of a complete graph.  Our construction
demonstrates also that in Wood's one-bend layout model, every graph can
be drawn in one or two layers, unless we introduce further restrictions
such as Wood's that the layout stay within a small area relative
to the vertex separation.

\section{Bounded geometric thickness}

\begin{figure}[t]
\centering
\includegraphics[width=3.5in]{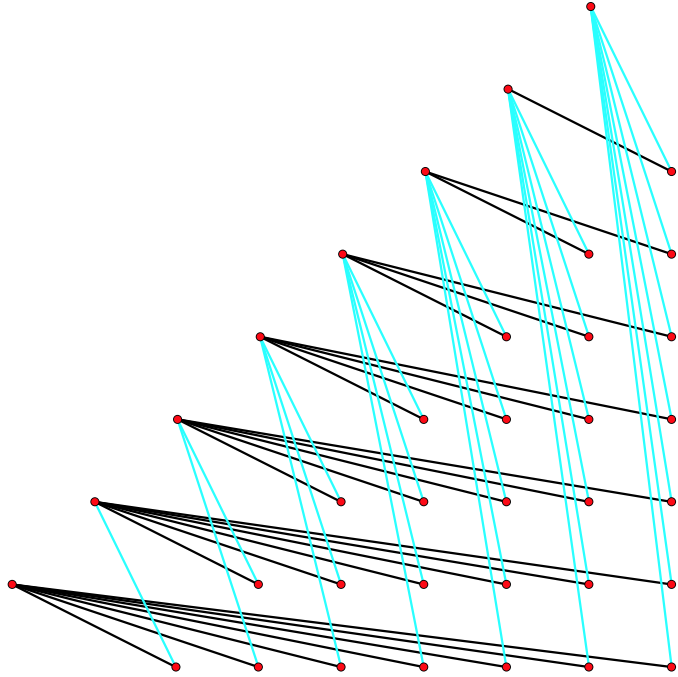}
\caption{Graph $G_8$ drawn with geometric thickness two.}
\label{fig:k8gt2}
\end{figure}

We define the graph $G_k$ by subdividing every edge of the complete
$k$-vertex graph $K_k$ into a path of two edges.
Equivalently, the vertices of $G_k$ can be viewed as corresponding to the
singleton and doubleton subsets of a $k$-element set, with an edge
between every two subsets having an inclusion relation.
Thus, $G_k$ has $k+{k\choose 2}$ vertices and $2{k\choose 2}$ edges.
Figure~\ref{fig:k8gt2} depicts $G_8$.
As the figure hints, all graphs $G_k$ can be drawn with small geometric
thickness:

\begin{theorem}
\label{thm:gtGk}
For $k\ge 5$, $\geomthickness(G_k)=2$.
\end{theorem}

\begin{proof}
Since $K_k$ and its subdivisions are nonplanar for $k\ge 5$, it is clear
that $\geomthickness(G_k)\ge 2$, so it remains to demonstrate the
existence of a two-layer drawing.
We let $v_i$, $0\le i<k$, denote the vertices of the complete graph $K_k$
from which $G_k$ is formed, and place vertex $v_i$
at the point with coordinates $(i,i+1)$.
To place the two-edge path between $v_i$ and $v_j$, $i<j$,
we use an edge on the first layer from $v_i$ to $(j+1,i)$ and an edge on
the second layer from that point to $v_j$.
All edges within a given row of the first layer, or within a given column
of the second layer, have a common endpoint, so there can be no crossings
within either layer.
\end{proof}

\section{Unbounded book thickness}

To show that $G_k$ does not have bounded book thickness, we need a
standard result of Ramsey Theory~\cite{GraRotSpe-80}, which we state in
the form we need:

\begin{lemma}
\label{lem:Ram}
For every pair of positive integers $c$ and $\ell$ there is an integer
$R_c(\ell)$ with the following property:
If the edges of the complete graph $K_{R_c(\ell)}$ are partitioned into
$c$ graphs, then at least one of the graphs contains a complete
$\ell$-vertex subgraph.
\end{lemma}

\begin{theorem}
For every positive integer $t$ there exists a $k$ such that
$\bookthickness(G_k)\ge t$.
\end{theorem}

\begin{proof}
Let $k=R_c(5)$, where $c={t-1\choose 2}$.  We then show that
$\bookthickness(G_k)\ge t$.
Suppose to the contrary that $G_k$ has a book embedding with
$t-1$ layers.  We use this drawing to partition the edges of $K_k$ into
$c$ subgraphs: an edge $v_iv_j$ is assigned to a subgraph
according to the unordered pair of layers used by the path connecting
$v_i$ to $v_j$ in
$G_k$.  If any of the two-edge paths in $G_k$ uses only a single layer,
the corresponding edge of $K_k$ may be placed arbitrarily into any of the
$t-2$ subgraphs involving that layer.

By Lemma~\ref{lem:Ram}, we can find a copy of $K_5$ in one of the
subgraphs.  This copy corresponds to a graph $G_5$, drawn with a book
embedding of only two layers.  However, $G_5$ is nonplanar, so its book
thickness is at least three, a contradiction.
\end{proof}

\section{Discussion}

We have shown that, for certain graphs $G_k$, the book thickness grows
arbitrarily while the geometric thickness is bounded.
It remains open how strongly separated these quantities are.
Our proof uses Ramsey theory, so yields only weak lower bounds on the
book thickness of $G_k$. The best upper bound we have been able to
achieve is $\bookthickness(G_k)=O(\sqrt k)$, by partitioning the vertices
of $K_k$ into blocks, using one layer per vertex within each block to
connect the vertices to reordered transfer points adjacent to the block,
and using one layer per block to connect the transfer points to vertices
in the other blocks.  We conjecture that this upper bound is tight.

The question of whether geometric thickness and thickness are
asymptotically equivalent is very interesting, and remains open.

\bibliographystyle{abuser}
\raggedright
\bibliography{thickness}
\end{document}